\documentclass{amsart}
\usepackage[square,comma]{natbib}
\usepackage[english]{babel} 
\usepackage{exscale}   
\usepackage{amsmath}
\usepackage{amsfonts}
\usepackage{amssymb}
\usepackage{amsxtra}
\usepackage{pstricks}
\usepackage{stmaryrd}
\usepackage{xspace}
\usepackage{epsfig}
\usepackage{mathrsfs}

\numberwithin{equation}{section}
\newtheorem{Theorem}{Theorem}[section]
\newtheorem{Lemma}[Theorem]{Lemma}
\newtheorem{Corollary}[Theorem]{Corollary}
\newtheorem{Proposition}[Theorem]{Proposition}
\newtheorem{Definition}[Theorem]{Definition}
\newtheorem{Remark}[Theorem]{Remark}

\newtheorem{Notation}[Theorem]{Notation}

\newcommand{\Cset}{\mathbb{C}}

\newcommand{\Nset}{\mathbb{N}}

\newcommand{\Sset}{\mathbb{S}}

\newcommand{\cA}{\ensuremath{{\mathcal A}}\xspace}         
\newcommand{\cB}{\ensuremath{{\mathcal B}}\xspace}         
\newcommand{\cE}{\ensuremath{{\mathcal E}}\xspace}         
\newcommand{\cF}{\ensuremath{{\mathcal F}}\xspace}         
\newcommand{\cH}{\ensuremath{{\mathcal H}}\xspace}         
\newcommand{\cM}{\ensuremath{{\mathcal M}}\xspace}         
\newcommand{\cN}{\ensuremath{{\mathcal N}}\xspace}         

\newcommand{\ctM}{\ensuremath{{\widetilde{\mathcal M}}}\xspace} 

\newcommand{\scrM}{\ensuremath{{\mathscr{M}}}\xspace} %
\newcommand{\ii}{\mathbf{i}}

\newcommand{\tvarphi}{\ensuremath{{\widetilde{\mathsf{\varphi}}}}\xspace} %
\newcommand{\1}{\ensuremath{{\rm 1\kern-.25em l}}\xspace}  
\newcommand{\alg}{{\operatorname{alg}}} 
\newcommand{\id}{\operatorname{id}}                        

\newcommand{\Ker}{\ensuremath{\operatorname{Ker}}\xspace}    
\newcommand{\Ran}{\ensuremath{\operatorname{Ran}}\xspace}    
\newcommand{\tail}{\ensuremath{\mathrm{tail}}}                 

\newcommand{\Aut}[1]{\ensuremath{\operatorname{Aut}(#1)}\xspace}  

\newcommand{\set}[2]{\mathopen{\{}#1\mathop{|}#2\mathclose{\}}}

\newcommand{\sotlim}{\ensuremath{\textsc{sot-}\lim}\xspace}   
\newcommand{\sot}{\ensuremath{\textsc{sot}}\xspace}           

\begin{document}
\title[On Lehner's `free' De Finetti theorem]{On Lehner's `free' noncommutative\\ analogue of De Finetti's theorem}
\author[C. K\"ostler]{Claus K\"ostler}
\address{University of Illinois at Urbana-Champaign, Department of Mathematics,
Altgeld Hall, 1409 West Green Street, Urbana, 61801, USA}
\email{koestler@uiuc.edu}
\date{June 16, 2008}
\begin{abstract}
Inspired by Lehner's results on exchangeability systems in \cite{Lehn06a} we define 
`weak conditional freeness' and `conditional freeness' for stationary processes in an 
operator algebraic framework of noncommutative probability. We show that these two 
properties are equivalent and thus the process embeds into a von Neumann algebraic 
amalgamated free product over the fixed point algebra of the stationary process.   
\end{abstract}
\maketitle
Recently Lehner introduced `weak freeness' for exchangeability systems within a cumulant
approach to *-algebraic noncommutative probability. 
A main result in \cite{Lehn06a} is that an exchangeability system with weak freeness and 
certain other properties embeds into an amalgamated free product analogous to the classical 
De Finetti theorem. Here we investigate Lehner's approach from an operator algebraic point 
of view which is motivated by recent results on a noncommutative version of De Finetti's 
theorem \cite{Koes08aPP} and a certain `braided' extension of it \cite{GoKo08a}. For the
classical De Finetti theorem, we refer the reader to Kallenberg's recent 
monograph \cite{Kalle05a} on probabilistic symmetries and invariance principles. 

Since tail events in probability theory lead to conditioning which goes beyond amalgamation,
we define `conditional freeness' and `weak conditional freeness' as a slight generalization 
of Voiculescu's `amalgamated freeness' \cite{Voic85a,VDN92a} and Lehner's `weak freeness' 
\cite{Lehn06a}, respectively. We investigate them for stationary processes in an operator 
algebraic setting of noncommutative probability (see \cite{Koes08aPP, GoKo08a} for details). 
Our main results are, re-formulated in terms of infinite minimal random sequences with 
stationarity:
\begin{enumerate}
\item[$\diamond$]
`Weak conditional freeness' and `conditional freeness' are equivalent; and the conditioning is
with respect to the tail algebra of the random sequence (see Theorem \ref{thm:character-freeness}).
\item[$\diamond$]
`Weak freeness' and `amalgamated freeness' are equivalent under a certain condition; and 
the amalgamation is with respect to the tail algebra of the random sequence 
(see Theorem \ref{thm:weak-amalgamated}).
\item[$\diamond$]
Each of these four variations of Voiculescu's central notion of freeness \emph{implies} that the
random sequence canonically embeds into a certain von Neumann algebra amalgamated free product
and that it enjoys exchangeability (see Theorem \ref{thm:character-freeness}, Theorem 
\ref{thm:weak-amalgamated} and their corollaries).
\end{enumerate}
Our results hint at that, dropping the assumption of stationarity, a certain asymptotic version of 
weak conditional freeness seems perhaps already to be equivalent to conditional freeness; which would
of course imply the distributional symmetry of exchangeability. 

We summarize the content of this paper. 
Section \ref{section:preliminaries} introduces `conditional freeness' and `weak 
conditional freeness' for stationary processes; and it provides a fixed point characterization
theorem from \cite{Koes08aPP}.  
Section \ref{section:main-results} contains our first main result, 
Theorem \ref{thm:character-freeness} on the equivalence of weak conditional freeness and 
conditional freeness. This results rests on an application of the mean ergodic theorem, also
provided there. Finally, we relate our results in Section \ref{section:application} to those 
obtained for exchangeability systems in \cite{Lehn06a}. We will see 
that, up to regularity and modular conditions, an exangeability system yields a stationary 
process in our sense. This observation is the starting point for Theorem \ref{thm:weak-amalgamated}, 
our second main result on the equivalence of amalgamated freeness and weak freeness. 
\section{Preliminaries}\label{section:preliminaries}
We are interested in a probability space $(\cM,\varphi)$ consisting of a 
von Neumann algebra $\cM$ with separable predual and a faithful normal state $\varphi$ on $\cM$.
A von Neumann subalgebra $\cM_0 \subset \cM$ is said to be \emph{$\varphi$-conditioned} if the $\varphi$-preserving conditional expectation $E_0 \colon \cM \to \cM_0$ exists. We say that
an endomorphism $\alpha \colon \cM \to \cM $ is \emph{$\varphi$-conditioned} if $\alpha$ is unital, $\varphi$-preserving 
and commutes with the modular automorphism group associated to $(\cM,\varphi)$. Throughout we will work in the GNS representation of $(\cM,\varphi$). Finally, $\bigvee_{i \in I}\cA_i$ denotes the von Neumann algebra generated by the family 
$(\cA_i)_{i \in I} \subset \cM$.  

\begin{Definition}\normalfont \label{def:rs}
A \emph{stationary process} $\scrM\equiv(\cM,\varphi, \alpha; \cM_0)$ consists of a probability space 
$(\cM,\varphi)$ which is equipped with a $\varphi$-conditioned endomorphism $\alpha$ and a $\varphi$-conditioned
von Neumann subalgebra $\cM_0 \subset \cM$. The \emph{canonical filtration} $\cF(\scrM)$ of $\scrM$ is the family 
$(\cM_I)_{I \subset \Nset_0}$ of subalgebras $\cM_I := \bigvee_{i \in I} \alpha^i(\cM_0)$.
We say that $\scrM$ is  \emph{minimal} if $\cM = \cM_{\Nset_0}$. 
\end{Definition}

\begin{Notation}\normalfont
The \emph{fixed point algebra} of $\alpha$ is denoted by $\cM^\alpha$
and $E$ is the $\varphi$-preserving conditional expectation from $\cM$ onto $\cM^\alpha$.
The \emph{tail algebra} of $\scrM$ is $\cM^\tail := \bigcap_{n \in \Nset_0} \alpha^n(\cM)$.
\end{Notation} 

\begin{Remark}\normalfont\label{rem:existence}
The $\varphi$-conditioning of $\cM_0$ and $\alpha$ imply that $\cF(\scrM)$ is family of 
$\varphi$-conditioned von Neumann subalgebras. In particular, the $\varphi$-conditioning of $\alpha$ 
ensures the existence of the $\varphi$-preserving conditional expectation $E$ from $\cM$ onto $\cM^\alpha$. 
We will make use of this in the proof of 
Theorem \ref{thm:character-freeness}.    
\end{Remark}

Motivated by Lehner's notion of `weak freeness' and the author's work on a 
noncommutative extended De Finetti theorem \cite{Koes08aPP,GoKo08a}, we 
introduce `weak conditional freeness' and `conditional freeness'. 

\begin{Definition}\normalfont \label{def:freeness}
Suppose $\scrM$ is a minimal stationary process.  
\begin{enumerate}
\item $\scrM$ and its filtration $\cF(\scrM)$ satisfy 
\emph{conditional freeness} if, for every $n \in \Nset$ 
and $n$-tuple $\ii\colon \{1,2,\ldots,n\} \to \Nset_0$,   
\[
E(x_1 x_2 x_3 \cdots x_n) = 0
\]
whenever 
\[
x_j \in (\cM^\alpha \vee \cM_{I_{\ii(j)}}) \,\cap\, \Ker E
\] 
with mutually disjoint subsets $\set{I_i}{i  \in \Ran \ii}$ and 
$\ii(1)\neq \ii(2) \neq \cdots \neq \ii(n)$.   
\item 
$\scrM$ and its filtration $\cF(\scrM)$ satisfy \emph{weak conditional freeness} 
if, for every $n \in \Nset$ and $n$-tuple $\ii\colon  \{1,2,\ldots,n\} \to \Nset_0$,  
\[
\varphi(x_1 x_2 x_3 \cdots x_n) = 0
\]
whenever
\[
x_j \in \cM^\alpha \vee \cM_{I_{\ii(j)}} \quad \text{and} \quad 
\lim_{N \to \infty}\frac{1}{N}\sum_{k=0}^{N-1}\varphi\big(x_j^* \alpha^k(x_j^{}) \big)=0
\]
with mutually disjoint subsets 
$\set{I_i}{i  \in \Ran \ii}$ and $\ii(1)\neq \ii(2) \neq \cdots \neq \ii(n)$.     
\end{enumerate}
\end{Definition}
\begin{Remark}\normalfont
Conditional freeness of $\cF(\scrM)$ is equivalent to Voiculescu's amalgamated freeness 
of the family $\big(\alpha^i(\cM_0\vee\cM^\alpha)\big)_{i \ge 0}$ in $(\cM,E)$. Thus it 
generalizes amalgamated freeness of two operators  \cite[Definition 3.8.2]{VDN92a} to 
amalgamated freeness of two families of operators. On the other hand, conditional freeness 
is a special case of conditional independence in \cite{Koes08aPP}. Note that the requirement 
$\cM^\alpha \subset \cM_0$ is very restrictive since $\cM_I \cap \cM^\alpha$ may be 
trivial for any finite set $I$. Such a situation occurs frequently for results of De Finetti 
type or, more generally, if tail events of random sequences are considered.   
\end{Remark}
\begin{Remark}\normalfont
Our notion of `conditional freeness' should not be confused with that given in \cite{BLS96a},
a generalization of free product states in a *-algebraic setting. 
\end{Remark}
\begin{Remark}\label{rem:weakfreeness}\normalfont
`Weak conditional freeness' formally simplifies to `weak freeness' 
(see \cite[Definition 4.1]{Lehn06a}) if
$\varphi\big(x_j^* \alpha^k(x_j^{}) \big) = \varphi\big(x_j^* \alpha^{N_0}(x_j^{}) \big)$ 
whenever $k>N_0$ for some $N_0 \in \Nset$. But there is also a significant difference: 
weak conditional freeness is formulated with respect to a family of von Neumann subalgebras 
\emph{saturated} by the fixed point algebra $\cM^\alpha$. In a *-algebraic or C*-algebraic 
approach one can not expect that a non-trivial fixed point of $\alpha$ is in $\cM^\alg := 
\alg\set{\alpha^n(\cM_0)}{n \in \Nset_0}$ or its norm closure; the situation 
$\cM^\alpha \cap \overline{\cM^\alg}^{\| \cdot \|} = \Cset \cdot \1_{\cM}$ may occur.    
\end{Remark}
We will need in Section \ref{section:application} a fixed point characterization result 
from \cite{Koes08aPP}.

\begin{Definition}\normalfont\label{thm:order-fact}
A stationary process $\scrM$ is said to be \emph{order $\cN$-factorizable} if 
$\cN$ is a $\varphi$-conditioned von Neumann subalgebra of $\cM$ and 
$E_{\cN}(xy) = E_{\cN}(x) E_{\cN}(y)$ for all 
$x\in \cM_I$ and $y \in \cM_J$ with $\max I < \min J $ or $\min I >\max J$. 
\end{Definition}
Note that the inclusion $\cN \subset \cM_I \cap \cM_J$ is \emph{not} required in this definition.  
\begin{Theorem}\label{thm:fixed-point}
Suppose the minimal stationary process $\scrM$ is order $\cN$-factorizable for 
the $\varphi$-conditioned von Neumann subalgebra $\cN$ of $\cM^\alpha$. Then it 
holds 
$$
\cN = \cM^\alpha =\cM^\tail.
$$  
\end{Theorem}
\section{Main result}\label{section:main-results}
We assume that the reader is familiar with von Neumann algebraic amalgamated 
free products. Their definition and the technical details of their construction 
can be found in \cite{JuPaXu07a} (see also \cite{VDN92a} for an outline).      
\begin{Theorem} \label{thm:character-freeness}
The following are equivalent for a minimal stationary process $\scrM$: 
\begin{enumerate}
\item[(a)]
$\cF(\scrM)$ satisfies weak conditional freeness;
\item[(b)]
$\cF(\scrM)$ satisfies conditional freeness;
\item[(c)]
$\cF(\scrM)$ embeds canonically into the von Neumann algebra amalgamated free product 
\[
(\ctM,\tvarphi):={\underset{\cM^\alpha}{*}}_{n=0}^{\infty}\big(\cM_0 \vee \cM^\alpha,\varphi|_{\cM_0 \vee \cM^\alpha}\big),
\]
such that the endomorphism $\alpha$ of $\cM$ is turned into the unilateral shift 
$\widetilde{\alpha}$ on the amalgamated free product factors of $\ctM$. 
\end{enumerate}
\end{Theorem}
We record an immediate consequence before giving the proof. Let $\Sset_\infty$ denote the
inductive limit of the symmetric groups $\Sset_n$. 

\begin{Definition}\normalfont
A stationary process
is said to be \emph{exchangeable} if, for any $n \in \Nset_0$, 
\[
\varphi\big(\alpha^{i_1}(a_1) \alpha^{i_2}(a_2)\cdots \alpha^{i_n}(a_n)\big)
=  \varphi\big(\alpha^{\pi(i_1)}(a_1) \alpha^{\pi(i_2)}(a_2)\cdots \alpha^{\pi(i_n)}(a_n)\big) 
\]
for all $n$-tuples $(i_1, \ldots i_n) \subset \Nset_0^n$ and $(a_1, \ldots a_n) \in \cM_0^n$ and
(finite) permutations $\pi \in \Sset_\infty$ on $\Nset_0$. 
\end{Definition}

\begin{Corollary}\label{cor:exchange}
A minimal stationary process $\scrM$ with weak conditional freeness is exchangeable.
\end{Corollary}
\begin{proof}
Due to Theorem \ref{thm:character-freeness} we can assume that $\scrM$ is already realized on
the von Neumann algebra amalgamated free product over the fixed point algebra of $\alpha$.
Now the transposition of the $(n-1)$-th and $n$-th factor in the amalgamated free product implements 
an $\varphi$-preserving automorphism $\gamma_n$ of $\cM$ and it is straight forward to verify that the 
subgroup $(\gamma_n)_{n \in \Nset}$ in \Aut{\cM} enjoys the relations
\begin{align*}
&&\gamma_i \gamma_{j} \gamma_i &= \gamma_{j} \gamma_i \gamma_{j} &\text{if $ \; \mid i-j \mid\, = 1 $,}&&  \\
&&\gamma_i \gamma_j &= \gamma_j \gamma_i  &\text{if $ \; \mid i-j \mid\, > 1 $,}&&  \\
&&\gamma_i^2 &= \id  &\text{for all $i \in \Nset$.} &&
\end{align*}
Thus we have found a representation of the symmetric group  $\Sset_\infty$ in $\Aut{\cM}$. 
The exchangeability of $\scrM$ is concluded from conditional freeness, $\varphi \circ \gamma_i = \varphi$ and 
$$
\alpha(x) = \sotlim_{n \to \infty} \gamma_1 \gamma_2 \cdots \gamma_n(x), \qquad x \in \cM,
$$
where $\sot$ denotes the strong operator topology. We conclude the exchangeability from 
\[
\varphi\big(\alpha^{i_1}(a_1) \alpha^{i_2}(a_2)\cdots \alpha^{i_n}(a_n)\big)
=  \varphi\big(\alpha^{\gamma_k(i_1)}(a_1) \alpha^{\gamma_k(i_2)}(a_2)\cdots \alpha^{\gamma_k(i_n)}(a_n)\big) 
\]
for all $k \in \Nset$, $n$-tuples $(i_1, \ldots, i_n) \in \Nset_0^n$ and $(a_1, \ldots a_n) \in \cM_0^n$
(see \cite{GoKo08a} for further details). This completes the proof.     
\end{proof}
We prepare the proof of Theorem \ref{thm:character-freeness} with an operator algebraic version of the mean ergodic theorem. 
\begin{Theorem}\label{thm:mean-ergodic}
Let $\scrM$ be a stationary process. Then for each $x\in \cM$,
\[
\sotlim_{N\to \infty} \frac{1}{N}\sum_{k=0}^{N-1} \alpha^k(x) = E(x).   
\]
\end{Theorem}
\begin{proof}
The strong operator topology ($\sot$) and the $\varphi$-topology generated by $x \mapsto \varphi(x^*x)^{1/2}$, $x\in \cM$,
coincide on norm bounded sets in $\cM$. Thus this mean ergodic theorem is an immediate consequence 
of the usual mean ergodic theorem in Hilbert spaces (see \cite[Theorem 1.2]{Pete83a} for example). 
\end{proof}

\begin{Corollary}\label{cor:weak-factor-i}
Suppose $\scrM$ is a minimal stationary process. 
Then for $x\in \cM$,
\begin{eqnarray*}
\lim_{N \to \infty} \frac{1}{N}\sum_{k=0}^{N-1}\varphi(x^* \alpha^k(x)) = 0
\qquad \Longleftrightarrow\qquad  
x \in \ker E. 
\end{eqnarray*}
\end{Corollary}
\begin{proof}
This is immediate from Theorem \ref{thm:mean-ergodic}, the faithfulness of $\varphi$ and $E$, and  
\begin{eqnarray*}
\lim_{N \to \infty} \frac{1}{N}\sum_{k=0}^{N-1}\varphi\big(x^* \alpha^k(x)\big) 
=  \varphi\big(x^* E(x)\big)
=   \varphi(E(x^*) E(x)\big)= 0.  
\end{eqnarray*}
\end{proof}

\begin{proof}[Proof of Theorem \ref{thm:character-freeness}]
`(a) $\Rightarrow$ (b)':
Let the tuple $(x_1, x_2, \ldots, x_n)$ satisfy 
the assertions of Definition \ref{def:freeness}(ii). Our goal is to show that then $(ax_1, x_2, \ldots, x_n)$ 
also satisfies them for any $a \in \cM^\alpha$. 
Due to Lemma \ref{cor:weak-factor-i} it suffices to prove 
\[ 
\lim_{N\to \infty} \frac{1}{N}\sum_{k=0}^{N-1}\varphi\big(x_1^* \alpha^k(x_1^{})\big) = 0  
\quad \Longrightarrow \quad
\lim_{N\to \infty} \frac{1}{N}\sum_{k=0}^{N-1}\varphi\big(x_1^*a^* \alpha^k(ax_1^{})\big) =0.
\]
Indeed the mean ergodic theorem,  the Kadison-Schwarz inequality and properties 
of conditional expectations yield
\begin{eqnarray*}
\lim_{N\to \infty} \frac{1}{N}\sum_{k=0}^{N-1}\varphi\big(x_1^*a^* \alpha^k(ax_1^{})\big) 
&=&  \varphi\big(x_1^*a^* E(ax_1^{})\big) \\
&=&  \varphi\big(E(x_1^*)a^* a E_{\cM^\alpha}(x_1)\big)\\      
&\le & \|a\|^2  \varphi\big(E(x_1^*) E(x_1)\big)\\
&=&  \|a\|^2 \lim_{N\to \infty}\frac{1}{N}\sum_{k=0}^{N-1}\varphi\big(x_1^* \alpha^k(x_1^{})\big).
\end{eqnarray*}
By our initial assumption, $(x_1, \ldots,x_n)$ satisfies
$
\lim_{N\to \infty}\frac{1}{N}\sum_{k=0}^{N-1}\varphi\big(x_1^* \alpha^k(x_1^{})\big)=0.  
$
We conclude from above estimates that, for any $a\in \cM^\alpha$, 
\[
\lim_{N\to \infty}\frac{1}{N}\sum_{k=0}^{N-1}\varphi\big(x_1^* a^*\alpha^k(ax_1^{})\big)=0.  
\]
Altogether we have shown that the tuple $(ax_1, x_2,\ldots,x_n)$ satisfies
all assertions of Definition \ref{def:freeness}(ii) if the tuple 
$(x_1, x_2,\ldots,x_n)$ does so. Thus, by weak conditional freeness,
$
\varphi(x_1 x_2 x_3 \cdots x_n) = 0
$
implies 
$\varphi(ax_1 x_2 x_3 \cdots x_n) = 0$ for every $a \in \cM^\alpha$. But this entails $E(x_1 x_2 x_3 \cdots x_n)=0$ by routine arguments, since we are working on bounded sets in $\cM$. \\ 
We note for the proof of `(b) $\Rightarrow$ (c)' that conditional freeness of $\cF(\cM)$ implies the amalgamated freeness of the family $(\cM_{\{i\}}\vee \cM^\alpha)_{i \in \Nset_0}$ in $(\cM,E)$ in the category of C*-algebraic probability spaces.
Thus Voiculescu's construction of the amalgamated free product of C*-algebras applies \cite{Voic85a}. Since the definition of $\scrM$ ensures the existence of $\varphi$-preserving conditional expectations (see Remark \ref{rem:existence}), all assumptions of \cite[Section 1]{JuPaXu07a} are satisfied which are needed for the construction of the von Neumann amalgamated free product $(\ctM,\tvarphi)$. Finally, one verifies that $\alpha$ becomes the shift $\widetilde{\alpha}$ on the amalgamated free product factors of $\ctM$. 
\\
The remaining implication `(c) $\Rightarrow$ (a)' follows from the fact that the shift $\widetilde{\alpha}$ has the
von Neumann subalgebra ${\underset{\cM^\alpha}{*}}_{n=0}^{\infty} \cM^\alpha$ of $\ctM$ as fixed point algebra. Now another
application of the mean ergodic theorem in the style of Corollary \ref{cor:weak-factor-i} shows the weak conditional freeness
of the minimal stationary process $\widetilde{\scrM}\equiv(\ctM,\tvarphi, \widetilde{\alpha}, \ctM_0)$, where $\ctM_0$ is the
canonical embedding of $\cM_0$.   
\end{proof}
\section{Application to Lehner's exchangeability systems} \label{section:application}
This section is devoted to the discussion under which conditions an exchangeability system (defined in \cite{Lehn06a}) leads to a stationary process $\scrM$ in the sense of Definition \ref{def:rs}. We will show that, under some mild regularity and modular conditions, Lehner's notion of weak freeness is equivalent to amalgamated freeness (over the fixed point algebra of the induced stationary process).   

We assume that the reader is familiar with the definitions and results of \cite{Lehn06a}. 
\begin{Definition}\normalfont
A \emph{*-algebraic probability space} $(\cA^\alg,\varphi^\alg)$ consists of a unital *-algebra over $\Cset$ and 
a unital positive linear functional $\varphi^\alg \colon \cA^\alg \to \Cset$. We say that $(\cA^\alg,\varphi^\alg)$ is a 
\emph{regular probability space} if:
\begin{center}
For every $a \in \cA^\alg$ there exists a constant $C_a$ such that $$\varphi^\alg(x^*a^*a x) < C_a \varphi^\alg(x^*x).$$ 
\end{center}
\end{Definition}
Note that this condition is satisfied if $\cA^\alg$ is a pre-C*-algebra. It is well known that such a regularity condition allows to extend the procedure of the GNS construction from C*-algebras 
to *-algebras over $\Cset$. Let us make this more precise for the convenience of the reader. 
\begin{Proposition}
Let $(\cA^\alg,\varphi^\alg)$ be a regular probability space. Then there exists a representation $\Pi$ of $\cA^\alg$ on
a Hilbert space $\cH$ and a vector $\xi \in \cH$ which is cyclic for $\Pi(\cH)$ such that 
$\varphi^\alg(x) = \langle \xi, \Pi(x) \xi \rangle$ for all $x \in \cA^\alg$.  
\end{Proposition}
\begin{proof}
Let $\cN:= \set{x \in \cA^\alg}{\varphi^\alg(x^*x)=0}$. The regularity condition implies that the left multiplication
$\cN \ni x \mapsto ax$ is bounded on $\cN$ for each $a\in \cA^\alg$. Thus the proof of the usual GNS construction literally translates. 
\end{proof}
We are ready to construct a stationary process $\scrM$ starting from a *-algebraic setting of infinite random sequences:  Suppose the regular probability space $(\cA^\alg, \varphi^\alg)$, a unital *-algebra $\cA_0^\alg$ 
(over $\Cset$) and  the family 
\[
(\iota_n^\alg)_{n \ge 0}\colon \cA_0^\alg \to \cA^\alg
\]
of unital *-algebra homomorphisms are given. 

Let $(\Pi, \cH, \xi)$ denote the GNS triple of the regular probability space $(\cA^\alg,\varphi^\alg)$. 
Then the double commutant $\cM := \Pi(\cA^\alg)^{\prime \prime}$ in $\cB(\cH)$ and $\varphi(x):= \langle \xi, x \xi \rangle$ define the probability space $(\cM,\varphi)$. 
We let $\cM_I:=  \bigvee_{i \in I} \Pi \circ\iota_i^\alg(\cA_0^\alg)$ and write $\cM_n$ if $I = \{n\}$. This gives us the filtration $(\cM_I)_{I\subset \Nset_0}$. Furthermore we may assume without loss of generality the minimality of this filtration. 

Now suppose that the family $(\iota_n^\alg)_{n \ge 0}$ is stationary, i.e. for every $N >0$,
$$
\varphi^\alg\big(\iota_{\ii(1)}^\alg(a_1)\cdots \iota_{\ii(p)}^\alg(a_p)\big) =  \varphi^\alg\big(\iota_{\ii(1)+N}^\alg(a_1)\cdots \iota_{\ii(p)+N}^\alg(a_p)\big) 
$$
for all $p$-tuples $\ii\colon \{1, \ldots, p\} \to \Nset_0$ and $(a_1, \ldots, a_p) \in \Big(\cA_0^\alg\Big)^p$. 
It follows from the asserted stationarity and minimality that there exists a unital $\varphi$-preserving endomorphism $\alpha$ of $\cM$ such that $\alpha^n(\cM_0) = \cM_n$. Consequently, the quadruple $(\cM,\varphi, \alpha; \cM_0)$ defines a (minimal) stationary process $\scrM$, provided two additional \emph{modular conditions} are satisfied:
\begin{enumerate}
\item[1$^\circ$]
$\cM_0$ is  $\varphi$-conditioned;  
\item[2$^\circ$]
$\alpha$ is  $\varphi$-conditioned.
\end{enumerate}   
We summarize above discussion. 
\begin{Lemma}\label{lem:discussion}
An exchangeability system $\cE$ in \cite{Lehn06a} yields a stationary process $\scrM$ if $\cE$ can be realized on a
regular probability space $(\cA^\alg, \varphi^\alg)$ such that the modular conditions 1$^\circ$ and 2$^\circ$ are
satisfied.  
\end{Lemma}
\begin{proof}
It is immediate from the definition of exchangeability (see \cite{Lehn06a} or \cite{Koes08aPP}) that exchangeability 
implies stationarity. The remaining construction follows above discussion. 
\end{proof}  

Lemma \ref{lem:discussion} motivates us to consider from now on a minimal stationary process $\scrM = (\cM,\varphi, \alpha; \cM_0)$ as the starting point for the further discussion of Lehner's weak freeness condition.  

\begin{Definition}[\cite{Lehn06a}]\normalfont \label{def:lehner}
Suppose $\scrM$ is a minimal stationary process with a weak*-dense *-algebra $\cM_0^\alg \subset \cM_0$ and let 
$\cM_I^\alg := \alg\set{\alpha^i(\cM_0^\alg)}{i \in I}$. We say that $\scrM$ and $\cF(\scrM)$ satisfy 
\emph{(algebraic) weak freeness}
if 
\[
\varphi(x_1 x_2 \cdots x_n) = 0
\]
whenever $\varphi(x_j^* \alpha^{N_j}(x_j))= 0$ for $x_j \in \cM^\alg_{I_j}$ and 
$N_j > \min \set{N}{x_j \in \cM^\alg_{\{0, \ldots,N\}}}$ with mutually disjoint subsets $\set{I_i}{i \in \Ran \ii}$
and $\ii(1) \neq \ii(2) \neq \cdots \neq \ii(n)$.     
\end{Definition}
\begin{Remark}\normalfont
The condition $N_j > \min \set{N}{x_j \in \cM^\alg_{\{0,\ldots, N\}}}$ simplifies to $N_j > \max I_j$ if the index set 
$I_j$ is bounded. Note that, for unbounded $I_j$, the choice of $N_j$ depends on $x_j$.  
\end{Remark}

We generalize this observation to our second main result which improves the main results in \cite{Lehn06a} for a large
class of exchangeability systems.  
\begin{Theorem} \label{thm:weak-amalgamated}
Let $\scrM$ be a minimal stationary process such that 
\begin{align}
\cM^\alpha \cap \cM^\alg_0 = \cM^\alpha \cap \cM^\alg_{\Nset_0}\tag{$*$} \label{$*$}
\end{align} 
for some weak*-dense *-algebra $\cM_0^\alg$ in $\cM_0$. Then the following are equivalent:
\begin{enumerate}
\item[(a)] $\cF(\scrM)$ satisfies weak freeness;
\item[(b)] $\cF(\scrM)$ satisfies amalgamated freeness in $(\cM, E)$;     
\item[(c)] $\cF(\scrM)$ embeds canonically into the von Neumann algebra amalgamated free product 
\[
(\ctM,\tvarphi):={\underset{\cM^\alpha}{*}}_{n=0}^{\infty}\big(\cM_0,\varphi|_{\cM_0}\big),
\]
such that the endomorphism $\alpha$ of $\cM$ is turned into the unilateral shift $\widetilde{\alpha}$ on the amalgamated free product factors of $\ctM$. 
\end{enumerate}
\end{Theorem}
An immediate consequence is that the assumption of `exchangeability' in \cite{Lehn06a} is turned into a conclusion.
\begin{Corollary}
A minimal stationary process with weak freeness is exchangeable. 
\end{Corollary}
\begin{proof}
Repeat the proof of Corollary \ref{cor:exchange}.
\end{proof}

\begin{proof}[Proof of Theorem \ref{thm:weak-amalgamated}]
`(a) $\Rightarrow$ (b)': Suppose $\scrM$ satisfies weak freeness. We conclude for $x_j$ (as stated in Definition \ref{def:lehner}) that
\begin{eqnarray*}
\varphi(x_j^* \alpha^{N}(x_j))
&=& \varphi(x_j^* \alpha^{N+1}(x_j))  
= \lim_{N \to \infty} \frac{1}{N}\sum_{n=0}^{N-1}\varphi(x_j^* \alpha^{n} (x_j))  \\  
\end{eqnarray*}
and, by Corollary \ref{cor:weak-factor-i},
\[
\varphi(x_j^* \alpha^{N}(x_j))= 0 \quad \Longleftrightarrow \quad E(x_j) =0. 
\] 
We identify in a second step which elements $x \in \cM^\alg_{\Nset_0}$ satisfy $E(x)=0$.
For this purpose let $\cN$ and $\cN_0$ be the von Neumann algebras generated by the orbit of 
$$
\cN^\alg:=\cM^\alpha \cap \cM^\alg_{\Nset_0} \quad \text{resp.} \qquad \cN_0^\alg:=\cM^\alpha \cap \cM^\alg_{0}
$$ 
under the action of the modular automorphism group associated
to $(\cM,\varphi)$. By construction and Takesaki's theorem, $\cN$ and $\cN_0$ are $\varphi$-conditioned
and we let $E_{\cN}$ resp.~ $E_{\cN_0}$ denote the corresponding conditional expectations 
from $\cM$ onto $\cN$ resp.~$\cN_0$.  Note that
$\cN \subset \cM^\alpha$ implies $E_{\cN} E= E_{\cN}$. Furthermore we know $\cN_0 \subset \cM_0$ 
(we do not know if $\cN \subset \cM_0$). 

Suppose $y \in \cM^\alg_{\Nset_0}$ and let $x := y - E(y)$. It is easy to see that $x \in \cM^\alg_{\Nset_0}$ 
if and only if $E(y) \in \cM^\alg_{\Nset_0}$. Thus $E(y) \in \cN^\alg \subset \cN$ and $E(y)= E_\cN E(y) = E_\cN(y)$. 
At this point we make use of the assumption \eqref{$*$} which ensures $\cN^\alg = \cN^\alg_0$, and consequently 
$\cN_0 = \cN$. Thus $x \in \cM^\alg_{\Nset_0}$ satisfies $E(x)=0$ if and only if 
$x = y -E_{\cN_0}(y)$ for some $y \in \cM^\alg_{\Nset_0}$.  

We conclude from this that the filtration $\cF(\scrM)$ satisfies amalgamated freeness in $(\cM, E_{\cN_0})$.
Suppose that we can prove that $\cN_0 \subset \cM^\alpha$ already implies $\cN_0 = \cM^\alpha$. Then Theorem \ref{thm:character-freeness} applies and we are done. But the implication that $\cN \subset \cM^\alpha$ forces 
$\cN_0 = \cM^\alpha$ is the content of the fixed point characterization result, Theorem \ref{thm:fixed-point}, 
as soon as we can ensure that amalgamated freeness in $(\cM, E_{\cN_0})$ implies order $\cN_0$-factorizability (see
Definition \ref{thm:order-fact}). This is easily verified and thus the proof is completed. 
\end{proof}

\begin{Remark}\normalfont
We do not know at the time of this writing if the assertion \eqref{$*$} 
is superfluous in Theorem \ref{thm:weak-amalgamated}. Can it be that
the *-algebra generated by $\cM^\alg_0$ and $\alpha(\cM^\alg_0)$ contains
more fixed points of $\alpha$ than $\cM^\alg_0$? 
\end{Remark}
\begin{Remark}\normalfont
The condition \eqref{$*$} can always be ensured by passing from the minimal stationary process
$\scrM$ to its saturation $\scrM_{\operatorname{sat}}:= (\cM, \varphi, \alpha; \cM_0 \vee \cN)$, where $\cN$ is as introduced in the proof of Theorem \ref{thm:weak-amalgamated}. Doing so the *-algebra
$\cM^\alg_0$ needs to be replaced by $\alg\{ \cM^\alg_0, \cN^\alg\}$. This procedure has the same effect as Lehner's
transfer from exchangeability systems to extended ones (see \cite[Remark and Definition 1.5]{Lehn06a}). 
\end{Remark}

\subsection*{Acknowledgments} The author thanks James Mingo and Roland Speicher for useful discussions on cumulants in free probability.

\bibliographystyle{alpha}                 
\label{section:bibliography}
\bibliography{claus26jan2008}                   


\end{document}